\documentclass{fpsac}
\usepackage{epsfig}

\newtheorem{theorem}{Theorem}[section]
\newtheorem{proposition}{Proposition}[section]
\newtheorem{lemma}[theorem]{Lemma}
\theoremstyle{definition}

\theoremstyle{remark}

\numberwithin{equation}{section}

\def\mT{\mathcal{T}}

\author[Sulanke]{Robert  A. Sulanke}
\address{Department of Mathematics, Boise State University,
    Boise, ID, USA 83725}

\email{sulanke@math.boisestate.edu}

\author[Xin]{Guoce Xin}
\address{Department of Mathematics, University of
Kentucky,
Lexington, KY, USA 40506}
\email{gxin@ms.uky.edu}
\title[Hankel Determinants for  Lattice Paths]{
 Hankel Determinants for Some Common Lattice Paths}

\subjclass[2000]{Primary 05A15; Secondary 15A36}

\keywords{Hankel determinants, lattice paths, continued fractions}

\begin{document}

\begin{abstract}
For a single value of $\ell$, let $f(n,\ell)$ denote the number of
lattice paths that use the steps $(1,1)$, $(1,-1)$, and
$(\ell,0)$, that run from $(0,0)$ to $(n,0)$, and that never run
below the horizontal axis.  Equivalently, $f(n,\ell)$ satisfies the
quadratic functional equation
  $F(x) = \sum_{n\ge 0}f(n,\ell) x^n
= 1+x^{\ell}F(x)+x^2F(x)^2.$
  Let $H_n$ denote the $n$ by $n$ Hankel
matrix, defined so that $[H_n]_{i,j} = f(i+j-2,\ell)$. Here we
 investigate the values of such determinants where
  $\ell = 0,1,2,3$.  For $\ell = 0,1,2$ we are able to
employ the Gessel-Viennot-Lindstr\"om method.  For the case $\ell=3$,
the sequence of determinants forms a sequence of period 14,
namely,
$$ (\det(H_n))_{n \ge 1} =
(1,1,0,0,-1,-1,-1,-1,-1,0,0,1,1,1,1,1,0,0,-1,-1,-1,\ldots)$$
For this case we are able to  use the continued
fractions method  recently introduced by Gessel and Xin.
 We also apply this
technique to evaluate Hankel determinants for other generating
 functions satisfying
a certain type of quadratic functional equation.

\begin{resume}
Pour une seule valeur de $\ell$,
soit  $f(n,\ell)$ le nombre des chemins treillis
que utilise les pas $(1,1)$, $(1,-1)$, et
$(\ell,0)$, vient de  $(0,0)$ \'a $(n,0)$, et que
ne vient jamais  dessous l'axis horizontale.
\'Equivalentement, le $f(n,\ell)$ satisfi\'e l'\'equation
fonctionnelle quadratique
  $F(x) = \sum_{n\ge 0}f(n,\ell) x^n
= 1+x^{\ell}F(x)+x^2F(x)^2.$
Soit $H_n$ le $n$ par $n$ matrice de Hankel,
d\'efinit pour que $[H_n]_{i,j} = f(i+j-2,\ell)$.
Nous examinons de tels d\'eterminants
o\'u
  $\ell = 0,1,2,3$.
Pour $\ell = 0,1,2$ nous pouvons employer
la m\'ethode de Gessel-Viennot-Lindstr\"om.
Pour le cas  $\ell=3$,
ls s\'equence de d\'eterminants forme une
s\'equence de p\'eriode $14$, \'a savoir
$$ (\det(H_n))_{n \ge 1} =
(1,1,0,0,-1,-1,-1,-1,-1,0,0,1,1,1,1,1,0,0,-1,-1,-1,\ldots)$$
Pour ce cas que nous pouvons utiliser
la m\'ethode de fractions continu\'ee
r\'ecemment introduit par Gessel et Xin.
Nous appliquons aussi cette technique
pour \'evaluer les d\'eterminants de Hankel
pour l'autres fonctions generatrices quie
satisfait un certain type d'\'equation fonctionnelle
qudratique.
\end{resume}

\end{abstract}

\maketitle  

\section{Introduction}

We will consider lattice paths  that use the
following three steps:
      $ U = (1,1)$, the up diagonal step;
     $ H = (\ell,0)$, the horizontal step of length $\ell$,
              where $ \ell$ is a single nonnegative
             integer;
  and $ D = (1,-1)$,   the down diagonal step.
Further, each $H$ step will be weighted by $t$,  and the others by
1. The weight of a path is the product of the weights of its
steps.  The weight of a path set is the sum of the weights of
its paths.

Let $f(n,t,\ell)$ denote the weight of the
path set   of
 paths  running from $(0,0)$ to $(n,0)$ that never run
below the $x$-axis.  When $t=1$, weight becomes cardinality.  For example,
\begin{itemize}
\item[-]   $f(n,0,0)$, equivalently  $f(n,0,\ell)$, is the weight of  a set of Dyck   paths, counted by the
aerated Catalan numbers:\\   $ (f(0,0,0), f(1,0,0),  f(2,0,0), \ldots  )  =
  ( 1, 0, 1,0, 2, 0, 5, 0, 14, 0, 42, 0, 132, 0, 429, \ldots ).
$

\item[-]    $f(n,t,1)$ is the weight of  a  set of Motzkin paths, counted by the Motzkin
numbers:\\  $ (f(0,1,1), f(1,1,1),  f(2,1,1), \ldots  )  =
 (1,  1, 2, 4, 9,
21, 51, 127, 323,835 \ldots ). $

\item[-]  $f(n,t,2)$ is the weight of  a  set
of  large Schr\"{o}der paths, counted by the
aerated large Schr\"{o}der numbers:\\ $ (f(0,1,2), f(1,1,2), f(2,1,2), \ldots
)  =
  ( 1, 0, 2, 0, 6, 0, 22, 0, 90, 0, 394, 0,  ).
$
\item[-]
$ (f(0,1,3), f(1,1,3),f( 2,1,3),\ldots  ) = (1, 0, 1, 1, 2, 3, 6, 10, 20,
36, 72, 136, 273, 532, \ldots )$

\end{itemize}
Previously, Pergola, {\it et al} \cite{Pergolaetal} and Sulanke
\cite{sulMomentsgenMot} have  considered such generalized Motzkin paths for
various values of
$\ell$ and have given additional references.
Letting $$F(x) = \sum_{n \geq 0} f(n,t,\ell) x^n$$  denote the generating
function
 for $f(n,t,\ell)$, we find by a common combinatorial decomposition that $F(x)$ satisfies the
 functional equation
$$ F(x) = 1 + t x^\ell F(x)  + x^2 F(x)^2.  $$

            Any sequence
 $A = (a_0, a_1, a_2 \ldots)$
 defines a sequence of  Hankel matrices, $ H_1, H_2,H_3
\ldots, $
 where  $H_n$ is an $n$ by $n$  matrix with entries
$(H_n)_{i,j} = a_{i+j-2}$.
For instance, the sequence $(f(n,1,3))_{n\ge0}$
yields
$$
H_1 = \left[\begin{matrix} 1 \end{matrix}\right],\quad
H_2 = \left[\begin{matrix} 1& 0\\ 0 & 1 \end{matrix}\right],\quad
H_3 = \left[\begin{matrix} 1& 0& 1\\ 0 & 1& 1\\  1& 1&2
\end{matrix}\right], \quad
H_4 = \left[\begin{matrix}
1& 0& 1& 1\\
                           0 & 1& 1&2 \\
                            1& 1&2&3\\
                            1& 2&3&6
                \end{matrix}\right]
$$

Our interest is to consider, for  nonnegative integer $\ell$,  the
corresponding sequence
of determinants $\det(H_n) $ where  each matrix $H_n$ has entries
$$(H_n)_{i,j} = f(i+j-2, t, \ell).$$
The following propositions constitute our main results:

\begin{proposition}\label{propmotzkin}
For $n \geq 0$,  $\ell = 1$, and arbitrary $t$ (including $t=0$,
yielding the Dyck path case)

          $$\det(H_n) =  1.$$
\end{proposition}

\begin{proposition}\label{propschr}
For $n \geq 0$,  $\ell = 2$, and arbitrary $t$ (including $t=0$,
yielding the Dyck path case),

 $$
\det(H_n) = \left\{
  \begin{tabular}{ll}
     $(1+t)^{n^2/4}$ & \mbox{ if n is even}\\
     $(1+t)^{(n-1)(n+1)/4}$ & \mbox{ if n is odd}
       \end{tabular}
\right.
$$
\end{proposition}

\begin{proposition}\label{propperiod14}
For  $t=1$ and $\ell = 3$,
$$ (\det(H_n))_{n \ge 1}^{14} = (1,1,0,0,-1,-1,-1,-1,-1,0,0,1,1,1 ).$$
Moreover, if $m,n \geq 0$ with $n-m = 0\mod 14$ then
 $\det(H_m) = \det(H_n)$.
\end{proposition}

In Section \ref{proofGVL}, using the  well-known  combinatorial
method of
 Gessel-Viennot-Lindstr\"om \cite{Bressoud}
 \cite{GV} \cite{Viennot}, we will prove
Propositions \ref{propmotzkin} and \ref{propschr}.  Our proof of
Propositions \ref{propmotzkin}  is essentially that of Viennot \cite{Viennot} who
also used the method to calculate various other Hankel
determinants relating to 
Motzkin paths.
Aigner \cite{Aigner} also studied such determinants. We note that
earlier Shapiro \cite{Shapiro} demonstrated that the Hankel
determinants for the usual Catalan numbers is 1. For the  large
Schr\"{o}der
 numbers $(r(n))_{n\ge0} = 1,2,6,22, 90, 394,
 \ldots$ whose generating function satisfies 
$$ R(x) = \sum_{k\ge 0} r(k)x^k = 1+xR(x)+xR(x)^2,$$
we show 
  that the  $n$-order Hankel   determinant is
$2^{n(n-1)/2}$, as stated in Proposition \ref{proplargesch}.

We remark that the problem of evaluating Hankel determinants
corresponding to a generating function has received significant
attention as considered  by Wall \cite{wall}. One of the basic
tools for such
evaluation 
 is the method of continued
fractions, either by $J$-fractions in Krattenthaler \cite{kratt} or  Wall
\cite{wall}  or by
$S$-fractions in Jones and Thron \cite[Theorem~7.2]{cfraction}.
However, both of these  methods need the condition that the determinant can
never be zero, a condition not always present in our study.
Recently, Brualdi and
Kirkland \cite{BrualdiKirkland}
used the $J$-fraction expansion to calculate   Hankel determinants for
various sequences related
to the
Schr\"{o}der numbers.  A slight modification of their proof of
\cite[Lemma 4.7]{BrualdiKirkland} proves our Proposition \ref{proplargesch}
for $t=1$.

 In Section \ref{proofperiod14} we
 establish the periodicity of 14 for the case $\ell = 3$ of
 Proposition \ref{propperiod14},
 by the continued fraction method  recently
 developed by Gessel and Xin \cite{xingessel2}.
 In the final
 section, we review  their technique more generally:  it
yields   a transformation for generating functions,
 satisfying a certain quadratic functional equation, that  also transforms
 the associated Hankel determinants in a simple manner.   We apply this
 transformation to evaluate the Hankel determinants for the cases
 $\ell = 1, 2$ (again) and for other path enumeration sequences related to
 $\ell = 3$.



\section{Employing the Gessel-Viennot-Lindstr\"om  method}\label{proofGVL}

Assuming a rudimentary knowledge of the Gessel-Viennot-Lindstr\"om  method,
we reformulate it to our needs.
All lattice paths use the three steps as previously  defined.
Given an $n$-tuple of lattice
paths on the $\mathbb{Z} \times \mathbb{Z}$ plane,
we say that  it is  {\it
nonintersecting} if no steps from different paths
share a common end point.  Thus an nonintersecting $n$-tuple may have paths 
crossing  or touching at points other than a common step end point.

Let $[(x_1,y_1), (x_2, y_2), \ldots, (x_n, y_n)]$ and
 $[(x'_1,y'_1), (x'_2, y'_2), \ldots, (x'_n, y'_n)]$
denote two lists of  distinct lattice points
such that
  $$x_{k+1} \leq x_{k} \le 0 \mbox{ and } 0 \le  y_{k}  \le  y_{k+1} $$
and
  $$0 \le x'_{k} \le x'_{k+1} \mbox{ and } 0 \le  y'_{k} \le   y'_{k+1}   .$$
\noindent We will refer to such  a  pair of lists as an ``{\sc
i-t-config}'' of order $n$ as their points will be the initial and terminal
points for each $n$-tuple of paths being   considered.

\begin{figure}[t]
     \begin{center}
     \leavevmode\psfig{figure=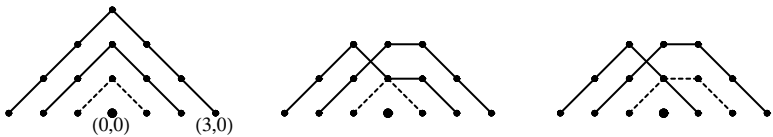,width=4in}
     \end{center}
\caption{\label{figure1}Some of the 4-tuples of paths for $\ell = 1$ and for
{\sc i-t-config} with $[ (0,0),(-1,0),(-2,0),(-3,0)]$ and
$[(0,0),(1,0),(2,0),(3,0)]$.  In each of these 4-tuples there is a
point path (a path of zero length) at $(0,0)$.  The first 4-tuple is
the only nonintersecting 4-tuple for this case.  The second and third
4-tuples are intersecting only at the point $(0,1)$. The second
4-tuple corresponds to the permutation $1243$
having sign of~-1, while the third corresponds to the permutation
$1342$
having sign of~1.  These two 4-tuples cancel one another under
the Gessel-Viennot-Lindstr\"{o}m  method. }
     \end{figure}

\begin{figure}[t]
     \begin{center}
     \leavevmode\psfig{figure=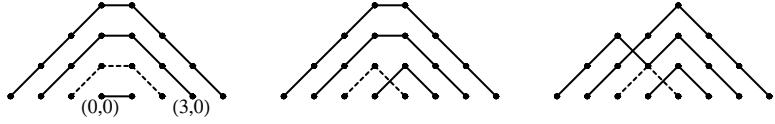,width=4in}
     \end{center}
\caption{\label{figure2}Three of the 4-tuples of paths for $\ell = 1$ and for
{\sc i-t-config} with $[ (0,0),(-1,0),(-2,0),(-3,0)]$ and
$[(1,0),(2,0),(3,0),(4,0)]$. The first and second 4-tuples are
both nonintersecting. The first has a signed weight of $t^4$ while the
second
has a signed weight of $-t^2$.  The third is intersecting only at the
point $(0,1)$.}
    \vspace{12pt}
     \end{figure}

 Let  $P_{i,j}$ denote the set of all  paths running  from
$(x_i,y_i)$ to $(x'_j,y'_j)$ that  never run   below the $x$-axis, with
  $|P_{i,j}|$ denoting the sum
 of the weights of its paths.
 Let $S_n$ denote the set of permutations
on $\{1, 2, 3, \ldots, n  \}$.  For any permutation $\sigma \in S_n$,
let $P_{\sigma}$ denote the set of all $n$-tuples of paths
$(p_1,p_2, \ldots, p_n)$, where $p_i \in P_{i,\sigma(i)} $ for
$1 \le i \le n$.
The {\it signed weight} of $(p_1,p_2, \ldots, p_n) \in P_{\sigma}$ is
 defined to be $\text{sgn}(\sigma)$ times the product of the weights of the
$n$ paths.  See Figures~\ref{figure1} and \ref{figure2}.

For our purpose the Gessel-Viennot-Lindstr\"om method is formulated in a
form similar to that in Viennot's 
notes \cite{Viennot}:
\begin{lemma}\label{Gessel-Viennot} Given an {\sc
i-t-config} of order $n$,
the sum of the signed weights of the nonintersecting $n$-tuples in
                     $\cup_{ \sigma \in  S_n} P_{\sigma}$
is equal to
                   $\det(\ (|P_{i,j}|)_{1 \leq i,j \leq n}\ )$.
\end{lemma}


\

\begin{proof}[
  Proof of Proposition \ref{propmotzkin}]  (A similar proof appears in
  \cite{Viennot}.) By Lemma \ref{Gessel-Viennot}
  $\det(H_n)$  is equal to the
sum of the signed weights of the nonintersecting $n$-tuples in
                     $\cup_{ \sigma \in  S_n} P_{\sigma}$
for the {\sc i-t-config} where $(x_i,y_i) = (-i+1,0)$ and
$(x'_i,y'_i) = (i-1,0)$, for $1 \leq  i \leq  n$.
Thus, for this {\sc i-t-config}, we seek  the  nonintersecting $n$-tuples.
First, the $1$-tuple $P_{1,1}$ contains just  the point
path beginning and ending at $(0,0)$. Next,  any
nonintersecting path from  $(-i+1,  0)$, for $1 < i \le n$,  must
begin with an $U$ step, while  any  nonintersecting path to
$(j-1,0)$, for $1 < j \le n$, must end with an $D$ step.
Repeating this analysis at
each integer-ordinate level $k$,  shows  the  nonintersecting path
from $(-i+1,0)$, $1\le i\le k$,
is forced to  be a sequence of $U$ steps followed by a sequence of
$D$ steps; moreover, it shows that any  nonintersecting path
from $(-i+1,0)$ to $(j-1,0)$, $k<i,j$,
must start with $k$  $U$ steps and end with  $k$  $D$ steps.
Inductively, each nonintersecting path is a sequence of $U$ steps followed
by a sequence of $D$ steps.
The $n$-tuple of such paths is the only
 nonintersecting $n$-tuple of $\cup_{ \sigma \in  S_n} P_{\sigma}$,
and it has weight equal $1$.
\end{proof}

\

We will use the following in  proving Proposition
\ref{propschr}:

\begin{lemma}\label{Lemma1}
For the lattice paths that use the steps $U$, $H=(2,0)$ , and $D$,
 that never run  below the x-axis, and
that have the {\sc i-t-config},
$$(x_i,y_i) = (-2i+2,0) \mbox{ and }
  (x'_i,y'_i) = (2i-2,0)
$$
for $1 \leq i \leq n$, the sum of the signed weights of the
nonintersecting $n$-tuples in $\cup_{\sigma \in S_n} P_{\sigma}$
equals  $(1+t)^{n(n-1)/2}$.
\end{lemma}

\begin{proof} For $(p_1,p_2,
\ldots, p_n) \in (P_{1,\sigma(1)} , P_{2,\sigma(2)}, \ldots,
P_{n,\sigma(n)} )$,
suppose  that   $(p_1,p_2, \ldots, p_n)$
is a nonintersecting $n$-tuple of paths for
some permutation $\sigma$.  Since the points in the {\sc
i-t-config} are spaced two units apart, the horizontal distance at
any integer ordinate between any two paths of $(P_{1,\sigma(1)} ,
P_{2,\sigma(2)},$ $ \ldots,$ $ P_{n,\sigma(n)} )$ must be even. It
follows inductively that, for $1 \leq i \leq n$, any path of
 the path set $P_{i,\sigma(i)}$ must begin with a sequence of
$i-1$ $U$-steps and finish with a sequence of $\sigma(i)-1$   $D$-steps. Thus,
computing the weight of the nonintersecting $n$-tuples is equivalent
to  computing the weight of the nonintersecting $n$-tuples for  the
new (``V'' shaped) initial-terminal configuration, denoted  by
{\sc i-t-config-new}, defined by
$$(x_i,y_i) = (-i+1,i-1) \mbox{ and }
  (x'_i,y'_i) = (i-1,i-1)
$$
for $1 \leq i \leq n$.

Before continuing, we notice,
for example when  $t=1$ and $n=4$, that the  matrix
$M(0)$ define by  $(M(0)_{i,j})_{1 \leq i,j \leq 4} =
(|P'_{i,j}|)_{1 \leq i,j \leq 4}$ for {\sc i-t-config-new}  is an
array of  Delannoy numbers. (See \cite{BanderierSchwer}, \cite{SulDel}.) When
$t=0$, $M(0)$ is the initial array from Pascal's triangle.
  In the following
array for $t=1$, the entries count the ways a chess king can move from the
north-west corner if it uses only
east, south, or south-east steps. Momentarily we will see
the role of the argument 0 in $M(0)$.
$$
M(0) = \left[
\begin{array}{cccc}
 1   & 1  & 1  & 1   \\
 1   & 3  & 5  & 7   \\
 1   & 5  & 13  &25      \\
 1   & 7  & 25  &63
 \end{array}
\right]
$$

Now for arbitrary $t$ and $n$,  let  $M(0)$ be the $n$ by $n$
matrix defined recursively by
$$ M(0)_{i,j} = M(0)_{i-1,j} +t  M(0)_{i-1,j-1} + M(0)_{i,j-1} $$
for  $1 < i$ and $1 <j$ with $ M(0)_{1,j} =1$ and $ M(0)_{i,1}
=1$ for  $1 \leq i$ and $1 \leq j$.   By Lemma \ref{Gessel-Viennot}
 $M(0) = (|P_{ij}|)_{1 \le i,j \le n}$ for  {\sc i-t-configNew}.
Thus
   $\det(M(0))$ is equal to the weight of the nonintersecting $n$-tuples for
  {\sc i-t-config}.  The proof is completed once we  show
$$\det(M(0)) = (1+t)^{n(n-1)/2}.$$

Given $ M(0)$,  we recursively define a sequence of $n$ by $n$
matrices
$$M(0),\  M(1 ),\  M(2 ),\ \ldots, M(n-1) $$
where, for $1 \leq k \leq n-1$,
$$
M(k)_{ij} = \left\{
\begin{tabular}{ll}
$M(k-1)_{i,j}$    &    \mbox{ for  $1  \leq  i  \leq k  $}  \\
$M(k-1)_{i,j} - M(k-1)_{i-1,j}$  &    \mbox{ for  $k+1  \leq  i
\leq n   $}
\end{tabular}
\right.
$$

With  {\sc claim}($k$) denoting the claim that
\begin{eqnarray*}
M(k)_{i,j} & = & M(k)_{i-1, j}+t  M(k)_{i-1, j-1}+  M(k)_{i, j-1}
\mbox{  for $i,j > k$,}\\
 M(k)_{i,i}& =&  (1+t)^{i-1}   \mbox{  for $i \leq k$,}\\
M(k)_{i,j } & =& 0   \mbox{  for $i > j$ and $j \leq k$,}\\
M(k)_{i,k+1} & =& (1+t)^{k}  \mbox{  for $ i \geq k+1 $,}
\end{eqnarray*}
one can establish {\sc claim}($k$) for $1 \leq k \leq n-1$ by
induction. Since  $M(n-1)$ is  upper triangular,
 we observe that
 $$\det(M(n-1)) = (1+t)^{n(n-1)/2}.$$
By  the type of  row operations used to obtain the sequence
 $M(0), M(1 ),  M(2 ), \ldots, M(n-1) $, their determinants are equal.
\end{proof}

Since, by the {\sc i-t-config} of
 Lemma \ref{Lemma1}, $(H)_{i,j} = |P_{i,j}|$ counts the large Schr\"{o}der paths from
 $(0,0)$ to  $ (2i+2j-2,0),$ immediately we have the
 the following corollary for the Hankel
determinants
of the weighted non-aerated
Schr\"{o}der numbers:
\begin{proposition}\label{proplargesch}
Let $f_n$ denote the weight of the path set of paths from $(0,0)$ to $(2n,0)$
which never run beneath the $x$-axis and where $H =(2,0)$ is weighted by $t$.
Equivalently, let $f_n$ satisfy
$$F(x) = \sum_{n\ge0}f_n x^n = 1 + txF(x) + xF(x)^2.$$
Then   the determinant of the  $n$-th order  Hankel matrix equals  $(1+t)^{n(n-1)/2}$.
\end{proposition}

 As a second  corollary to Lemma \ref{Lemma1}, we have
\begin{lemma}\label{Lemma2}
For the lattice paths that use the steps $U$, $H=(2,0)$, and $D$,
that never run  below the x-axis, and that have the {\sc
i-t-config} with
$$(x_i,y_i) = (-2i+1,0) \mbox{ and }
  (x'_i,y'_i) = (2i-1,0)
$$
for $1 \leq i \leq n$, the sum of the signed weights for the
nonintersecting $n$-tuples in $\cup_{\sigma \in S_n} P_{\sigma}$ is
$(1+t)^{n(n+1)/2}$.
\end{lemma}

\begin{proof}[Proof of Lemma \ref{Lemma2}] We first translate all paths upwards
one unit and then prepend
 a $U$-step and append a $D$-step to every path.  Next we
add the point path at $(0,0)$. The sum of the signed  weights of
the nonintersecting $n$-tuples in the original
 configuration equals
that of  the nonintersecting $n+1$-tuples in this new configuration,
which in turn is given by Lemma \ref{Lemma1}.
\end{proof}

\begin{proof}[Proof
of Proposition \ref{propschr}] Suppose that $n$ is even; the
proof when $n$ is odd is similar.
Here the Hankel matrix $(|P_{i,j}|)_{1\le i,j\le n}$ corresponds to the
 {\sc i-t-config} with
 $$(x_i,y_i)= (-i+1,0) \mbox{ and } (x'_i,y'_i)=  (i-1,0) \mbox{ for }
 1\le i \le n.$$
Since $\ell = 2$,  no endpoint of a
step on a path that originates from an oddly indexed
 initial point (i.e., a point $(-i+1,0)$ for odd $i$)
 will intersect an  endpoint of a step on a path that originates from
an evenly indexed initial point.  Moreover, for any permutation
$\sigma$ corresponding to a nonintersecting $n$-tuple,
$\sigma(i)-i$ must be even for each $i$, and hence $\text{sgn}(\sigma) =1$.
Thus the weight of the
nonintersecting $n$-tuples is the product of the weight of those
originating from oddly indexed initial points times the weight of
those originating from evenly indexed initial points.

Hence, with  $m = n/2$,
   let {\sc i-t-configA} have
 $$
(x_i,y_i) = (-2i+2,0) \mbox{ and }
  (x'_i,y'_i) = (2i-2,0)
\mbox{ for $1 \leq i \leq m$,}
$$
and let  {\sc i-t-configB}
have
$$
(x_i,y_i) = (-2i+1,0) \mbox{ and }
  (x'_i,y'_i) = (2i-1,0)
\mbox{ for $1 \leq i \leq m$.}
$$
Applying Lemmas \ref{Lemma1} and \ref{Lemma2} to these configurations
yields the weight of
nonintersecting $n$-tuples of the original configuration as
$$
(1+t)^{m(m-1)/2} (1+t)^{m(m+1)/2} = (1+t)^{m^2}  =
(1+t)^{n^2/4}.$$
\end{proof}

\

Next we consider Hankel determinants for sequences of path weights that
ignore the initial term. For the sequence  $ f(1,t,\ell),
f(2,t,\ell), \ldots, $ we will let $H^1_n$  denote the
matrix  where the entries satisfy
$(H^1_n)_{i,j} = f(i+j-1,
t, \ell)$. See  Figure~\ref{figure2}.

\begin{proposition}\label{L1missing}
For  $ \ell = 1$ (Motzkin case again), the sequence of
determinants satisfies the recurrence
$$  \det(H^1_n) = t\ \det(H^1_{n-1}) -  \det(H^1_{n-2 }) $$
subject to  $  \det(H^1_1) = t$ and $  \det(H^1_2)  =  (t-1)(t+1)$.
\end{proposition}

\begin{proof}    Aigner \cite{Aigner}
considered the case for $t=1$.  For arbitrary $t$, our
 proof considers how the particular paths must
look  in the nonintersecting case.
Observe that $\det(H^1_n)$ is the sum of the weights of
the nonintersecting $n$-tuples
for the {\sc
i-t-config}($n$) taken as
$$[(0,0),(-1,0),\ldots,(n-1,0)] \mbox{ and } [(1,0),(2,0),\ldots,(n,0)].$$
Each of these nonintersecting $n$-tuples
 belongs to one of two  types:
(1) those
containing the path from $(0,0)$ to $(1,0)$ with all other paths forced to
begin with $U$, end with $D$, and have ordinate at least one elsewhere;  (2)
those  containing the path $UD$ from $(0,0)$ to $(2,0)$
and the path  $UD$ $(-1,0)$ to $(1,0)$ with all other paths forced to
begin with $UU$, end with $DD$, and have ordinate at least two elsewhere.
The set  of the first type has a total weight $t$ times the  sum of the weights of
the nonintersecting $(n-1)$-tuples
on the {\sc
i-t-config}($n$-1), which is
$t \det(H^1_{n-1})$. Since each $n$-tuple of the second   type has
the defined crossing of the  path from  $(0,0)$ with  that from $(-1,0)$, the set has total weight is the sign of the
corresponding permutation  times the  sum of the weights of
the nonintersecting $(n-2)$-tuples
on the {\sc
i-t-config}($n$-2), which is
$-\det(H^1_{n-2})$.
\end{proof}

For $\ell = 2$, we will indicate how Lemma \ref{Lemma2}
proves

\begin{proposition}\label{L2missing}

For $n \geq 0$,  $\ell = 2$, and arbitrary
$t$,  the sequence of determinants satisfies
$$
 \det(H^1_n) = \left\{
  \begin{tabular}{ll}
     $0$   & \mbox{ if n is odd}\\
     $(-1)^{n/2}(1+t)^{n(n+2)/4}$ & \mbox{ if n is even}
  \end{tabular}
\right.
$$
\end{proposition}
\begin{proof}
Here the Hankel matrix  can correspond to
 {\sc i-t-config} with
 $$(x_i,y_i)= (-i+1,0) \mbox{ and } (x'_i,y'_i)=  (i,0) \mbox{ for }
 1\le i \le n.$$
Since $\ell = 2$, if there is a path from $(x_i,y_i)$ to $(x'_j,y'_j)$,
then $i-j$ is odd.  It follows that, if $n$ is odd, there can be no
$n$-tuples of paths for the configuration.  If $n$ is even and $m = n/2$,
the sign of any permutation for an nonintersecting  $n$-tuples can be shown to
be
$(-1)^m$.
Thus the weight of the
nonintersecting $n$-tuples is $(-1)^m$ times the weight of those
originating from oddly indexed initial points times the weight of
those originating from evenly indexed initial points.
The proof is completed by applying \ref{Lemma2} to
{\sc i-t-configA} with
 $$
(x_i,y_i) = (-2i+2,0) \mbox{ and }
  (x'_i,y'_i) = (2i,0)
\mbox{ for $1 \leq i \leq m$,}
$$
and     to {\sc i-t-configB}
with
$$
(x_i,y_i) = (-2i+1,0) \mbox{ and }
  (x'_i,y'_i) = (2i-1,0)
\mbox{ for $1 \leq i \leq m$.}
$$
\end{proof}

\section{Periodicity fourteen and continued fractions}\label{proofperiod14}

Here we will  repeated apply the ``continued fractions method''
recently developed by Gessel
and Xin
\cite{xingessel2} to determine the periodicity of
the sequence of Hankel determinants  for $\ell=3$ and $t=1$.
  This method,  presented more formally in the next section,
transforms both generating functions and corresponding determinants.
In this section we will concentrate on the specific generating function $F(x)$
satisfying
$$F(x)= 1+x^3 F(x) +x^2 F(x)^2.$$
{}From this functional equation, or from the related recurrence
for its coefficients,  there appears to be no clue why the
associated sequence of Hankel determinants should have a period of
14.

 For an arbitrary generating function $D(x,y)=\sum_{i,j=0}^\infty d_{i,j}
x^i y^j$, let $[D(x,y)]_n$ denote the $n$ by $n$
determinant
$\det((d_{i,j})_{0\le i,j\le n-1}).$
For any  $A(x)=\sum_{n\ge
0} a_n x^n$, define the \emph{Hankel matrix for} $A$ {\it of order} $n$, $n \ge 1$,
by $H_n(A)=(a_{i+j-2})_{1\le i,j\le n}$. It is straight forward to show that the
 Hankel determinant  $\det(H_n(A))$ can be expressed as
$$\det(H_n(A)) = \left[ \frac{xA(x)-yA(y)}{x-y}\right]_n.$$

We will use an easily-proven ``product rule'' of \cite{xingessel2}
for transforming the
generating functions:
 {\it If $u(x)$ is a formal power series with $u(0)=1$, then}
$$[u(x)D(x,y)]_n =[D(x,y)]_n=[u(y)D(x,y)]_n.$$
We will make five   transformations showing, for $n \ge 8$,
$$\det(H_n(F)) =
     \det \left( \text{diag}\left(  [1],[1], \left[ \begin {array}{ccc}
     0&0&1\\0&1&0
\\1&0&-2\end {array} \right],[1],[1],H_{n-7}(F)\right)\right),$$
where the right side is the determinant of a block-diagonal matrix
consisting of six blocks along the diagonal, four of which are 1 by 1
identity matrices, and having entry 0 elsewhere.
It then follows that $\det(H_n(F))=-\det(H_{n-7}(F))$. This
implies
    that the period for $\det(H_n(F))$ is 14, and Proposition
    \ref{propperiod14} will
    be proved.

We start with $F_0(x)=F(x)$, and define $F_i(x)$ from $F_{i-1}(x)$
according to a transformation where each  Hankel
determinant for $F_i(x)$ are derived from one for
$F_{i-1}(x)$ with the aid of the product rule, which is not always mentioned.
In the following,
 $F_i(x)$ will always satisfy a quadratic functional
equation
$$a(x)F_i(x)^2+b(x)F_i(x)+c(x)=0,$$
which  is equivalent to the continued fraction
$$F_i(x)= \frac{-c(x)}{b(x)+a(x)F_i(x)}. $$
In particular, for $\ell = 3$,  $$F_0(x)= \frac{1}{1-x^3-x^2F_0(x)}.$$

{\bf Transformation 1:} Using this continued fraction of $F_0$,
 substitution,
and simplification we obtain
$$  \det(H_n(F)) =
\left[ \frac{xF_0(x)-yF_0(y)}{x-y}  \right]_n  =
\left[ \frac{-x{y}^{2}F_0(y)+y{x}^{2}F_0(x)+ \left( x-y \right)
\left( y{x}^{2}+x{y}^{2}+1
 \right) }{(1-x^3-x^2F_0(x))(1-y^3-y^2F_0(y))(x-y)}  \right]_n.$$
Multiplying by $(1-x^3-x^2F_0(x))(1-y^3-y^2F_0(y))$, which will
not affect  the value of the determinant   by the product rule, we can write the
determinant as
$$ \left[ 1+xy \frac{xF_1(x)-yF_1(y)}{x-y } \right]_n$$
where
\begin{equation}\label{transform1} F_1(x)=F_0(x)+x. \end{equation}
The associated matrix is block-diagonal with two blocks: the matrix
$[1]$ and  the Hankel matrix for $F_1(x)$.
Certainly,
$$ \det(H_n(F_0))= \det(H_{n-1}(F_1)).$$
{}From (\ref{transform1}) and the functional equation for
$F_0(x)$, we obtain  the  functional equation
$$F_1(x)= {\frac {1+x}{1+{x}^{3}-{x}^{2}F_1(x)}}.$$

{\bf Transformation 2:}
Using this continued fraction for  $F_1$, substituting in
  $\displaystyle \frac{xF_1(x)-yF_1(y)}{x-y}$, and  multiplying by
$(1+x^3-x^2F_1(x))(1+y^3-y^2F_1(y)) $      yields
$$ \left[ \frac{xF_1(x)-yF_1(y)}{x-y} \right]_n= \left[ \frac{-x{y}^{2} \left( x+1 \right) F_1(y)+y{x}^{2} \left( y+1
  \right) F_1(x)
- \left( y +1 \right)  \left( x+1 \right)  \left( xy-1 \right)
\left( x-y
 \right)}{x-y}  \right]_n.$$
Upon multiplying by $(1+x)^{-1}(1+y)^{-1}$, the
determinant is equal to
$$ \left[  1+ xy\frac{xF_2(x)-yF_2(y)}{x-y} \right]_n,$$
where
\begin{equation}\label{transform2} F_2(x) =F_1(x)/(1+x)-1. \end{equation}
The associated matrix being block diagonal shows
$$\det (H_{n-1}(F_1)) = \det(H_{n-2}(F_2)).$$
 {}From (\ref{transform2})
and the functional equation for $F_1(x)$, we obtain
$$F_2(x) ={\frac {{x}^{2}}{1-2{x}^{2}-x^3- \left( {x}^{3}+{x}^{2}
 \right) F_2(x)}}.$$

{\bf Transformation 3:}
Substituting for $F_2$ with the above fraction,  simplifying,
and multiplying by $(1+x)(1-x-x^2-x^2F_2(x))(1+y)(1-y-y^2-y^2F_2(y))$ shows that the determinant
$ \displaystyle \left[
\frac{xF_2(x)-yF_2(y)}{x-y} \right]_n $ equals $$ \left[\frac{{y}^{2}{x}^{3} \left( y+1 \right) F_2(y)-{x}^{2}{y}^{3} \left( x+1 \right)
{}F_2(x)- \left( x-y \right)  \left(
2\,{y}^{2}{x}^{2}-{x}^{2}-xy-{y}^{2}
 \right)}{x-y}  \right]_n $$
which can be rewritten as
$$\left[  x^2+xy+y^2  -2x^2 y^2 +x^3y^3\frac{xF_3(x)-yF_3(y)}{x-y}  \right]_n ,
$$
where $F_3(x)$ is  indeed a power series satisfying
\begin{equation}\label{transform3}
F_3(x) =(x+1)F_2(x)/x^2.
\end{equation}
This time the corresponding matrix is a block-diagonal matrix with the
block $ \left[ \begin {array}{ccc} 0&0&1\\0&1&0
\\1&0&-2\end {array} \right]$ followed by the Hankel matrix for
$F_3(x)$. Hence
$$\det(H_{n-2}(F_2)) = -\det(H_{n-5}(F_3)). $$
{}From (\ref{transform3}) and the functional equation for
$F_2(x)$, we obtain $$ F_3(x)= {\frac
{1+x}{1-2x^2-{x}^{3}-{x}^{4}F_3(x)}}.$$

{\bf Transformation 4:} Substituting for $F_3$ with the fraction,  simplifying,
and multiplying by $(1-2x^2-x^3-x^4F_3(x) )((1-2y^2-y^3-y^4F_3(y) ) $  the
determinant
$\displaystyle  \left[ \frac{xF_3(x)-yF_3(y)}{x-y}\right]_n $ equals
$$\left[ \frac{-x{y}^{4} \left( x+1 \right) F_3(y)+y{x}^{4} \left( y+1 \right) F_3(x)+
  \left( y
+1 \right)  \left( x+1 \right)  \left( xy+1 \right)  \left( x-y
 \right) }{x-y}  \right]_n . $$
By multiplying the generating function  by $(1+x)^{-1}(1+y)^{-1}$, this determinant becomes
$$\left[ 1+xy \frac{xF_4(x)-yF_4(y)}{x-y}\right]_n ,$$
where
\begin{equation}\label{transform4} F_4(x) =1+x^2{}F_3(x)/(1+x).
 \end{equation}
 Therefore,
$$\det(H_{n-5}(F_3))=\det(H_{n-6}(F_4)).$$
{}From (\ref{transform4}) and the functional equation for
$F_3(x)$, we obtain
$$F_4(x)= \frac{1}{1+x^3 -\left( {x}^{3}+{x}^{2} \right) F_4(x) }.$$

{\bf Transformation 5:} Substituting for $F_4$ with the above
fraction,  simplifying, and multiplying by $(1-x^2F_4(x)
)(1-y^2F_4(y) ) $  the determinant $\displaystyle  \left[
\frac{xF_4(x)-yF_4(y)}{x-y} \right]_n $ equals
$$  \left[ \frac{-x{y}^{2} \left( y+1 \right) F_4(y)+{x}^{2}y \left( x+1
  \right) F_4(x)
- \left( x -y \right)  \left( y{x}^{2}+x{y}^{2}-1 \right) }{x-y}\right]_n
= \left[ 1+xy \frac{xF_5(x)-yF_5(y)}{x-y} \right]_n ,$$
where $F_5(x)=(1+x)F_4(x)-x$.  Hence, $\det(H_{n-6}(F_4))=\det(H_{n-7}(F_5)).$

Finally, it is routinely  checked that $F_5(x)=F_0(x)$.

\section{The quadratic transformation for Hankel determinants}
\def\mT{\mathcal{T}}

One can use the method introduced in the
previous section  to evaluate the
Hankel determinants for generating functions satisfying a certain type of
quadratic functional equation. The generating functions $F(x)$ in this section are the unique
solution of a quadratic functional equation satisfying \begin{equation} \label{e-F} F(x)= \frac{x^d}{u(x)+x^k
v(x)F(x)},
\end{equation}
where $u(x)$ and $v(x)$ are rational power series with nonzero
constants, $d$ is a nonnegative integer, and $k$ is a positive
integer. Note that if $k=0$, $F(x)$ is not unique.
Our task now is to
derive a transformation $\mT$ so that $\det( H_n(F)) =a\det( H_{n-d-1}(\mT (F)))$
for some  value $a$ and nonnegative integer $d$.
In addition to Hankel matrices for the power series $A=\sum_{n\ge 0}a_ix^i$,
 we will  consider  \emph{shifted
Hankel matrices}:  $H_n^k(A)$ denotes the matrix
$(a_{i+j+k-2})_{1\le i,j\le n}$.   Shifted matrices have
  appeared in Proposition \ref{L1missing} and \ref{L2missing}.

The first  proposition is elementary:
\begin{proposition}\label{l-constantTr}
If $F$ satisfies \eqref{e-F}, then $G=u(0)F$ satisfies
$$\det
(H_n(G))=u(0)^n \det (H_n(F)), \text{ and }\ \ G(x)
=\frac{x^d}{u(0)^{-1} u(x) +x^k u(0)^{-2} v(x)G(x)}.$$
\end{proposition}

\begin{proposition}\label{l-quadraticTr}
Suppose $F$ satisfies \eqref{e-F} with $u(0)=1$. We separate
$u(x)$ uniquely as $u(x)=u_L(x)+x^{d+2} u_H(x)$, where $u_L(x)$ is a
polynomial of degree at most $d+1$ and $u_H(x)$ is a power
series.
\begin{enumerate}

\item[(i)] If $k=1$, then there is a unique $G$ such that
$$ G(x) =\frac{-v(x)-xu_L(x)u_H(x)}{u_L(x)-x^{d+2}u_H(x)-x^{d+1}G(x)},$$
Moreover, $$G(x)=-xu_H(x)-x^{-d} v(x)F(x)$$
and a shifted matrices appears with
$$\det (H^1_{n-d-1}(G(x)))=(-1)^{d(d+1)/2}\det (H_n(F(x))).$$
\item[(ii)] If $k\ge 2$, then there is a unique $G$ such that
$$ G(x) =
  \frac{-x^{k-2}v(x)-u_L(x)u_H(x)}{u_L(x)-x^{d+2}u_H(x)-x^{d+2}G(x)},$$
Moreover, $$G(x)=-u_H(x)-x^{k-d-2}v(x)F(x) $$
and
$$ \det (H_{n-d-1}(G(x)))=(-1)^{d(d+1)/2}\det (H_n(F(x))). $$

\end{enumerate}
\end{proposition}

\begin{proof}
We prove only part (ii) as part (i) is similar.
The generating function for $H_n(F)$ is given by
\begin{align*}
\frac{xF(x)-yF(y)}{x-y}& =\frac{1}{x-y}\left(
\frac{x^{d+1}}{u(x)+x^k
  v(x)F(x)}-\frac{y^{d+1}}{u(y)+y^k v(y)F(y)}\right)\\
&= \rule{0pt}{18pt}
\frac{-y^{d+1}u(x)-y^{d+1}x^kv(x)F(x)+x^{d+1}u(y)
+x^{d+1}y^kv(y)F(y)}{\left(u(x)+x^k v(x)F(x)\right)
\left(u(y)+y^k v(y)F(y)\right) (x-y)}
\end{align*}

\noindent We can multiply by $(u(x)+x^kv(x)F(x))$ and by $(u(y)+y^kv(y)F(y))$
without changing the above determinant by the product rule. Next we
observe that $x^d$ divides $F(x)$, and write
$u(x)=u_L(x)+x^{d+2}u_H(x)$ as in the proposition. The  resulting generating
function can be written as

$$\frac{-y^{d+1}u_L(x)+x^{d+1}u_L(y)}{x-y}+ (xy)^{d+1}
\frac{-x(u_H(x)-x^{k-d-2}v(x)F(x))  +y(u_H(y)+y^{k-d-2}v(y)F(y))}{x-y}.$$

We now set $G(x)=-u_H(x)-x^{k-d-2}v(x)F(x) $, which can be straightforwardly
shown to agree with the  defining functional equations.
Suppose that
$u_L(x)=1+a_1x+\cdots a_{d+1}x^{d+1}$, then $\left[
\displaystyle{\frac{xF(x)-yF(y)}{x-y}}
 \right]_n$
is equal to the determinant of the block-diagonal matrix
$$\text{diag} \left(\left(\begin{array}{cccc}
0 & \cdots & 0 & 1 \\
\vdots&\vdots & \vdots & \vdots \\
0 & 1 &\vdots & a_{d-1} \\
1  & a_1 &\cdots & a_d
\end{array}\right), H(G(x)) \right).$$
\noindent The determinant of the first block is easily seen to be
$(-1)^{d(d+1)/2}$.
\end{proof}

\

Given these  propositions and that $H^1(A)=H(x^{-1}(A(x)-A(0)))$
for any series $A$,
we can
now define our transformation $\mT(F)$: {\it
For $F$ satisfying \eqref{e-F},
\begin{itemize}
\item if $u(0)\ne 1$, then $\mT(F)=G$, as given in Proposition
\ref{l-constantTr}.
\item if $u(0)=1$ and $k=1$, then $\mT(F)= x^{-1}(G(x)-G(0))$,
with $G$ given in Proposition \ref{l-quadraticTr}(i).
\item  if $u(0)=1$ and $k\ge 2$, then $\mT(F)=G$, as given in Proposition
\ref{l-quadraticTr}(ii).
\end{itemize}
Moreover, the relation between $\det( H_n(F))$ and $\det (H_n(\mT(F)))$
is given in Propositions \ref{l-constantTr} and \ref{l-quadraticTr}.}

\

{\bf Example 1:} {\sc Other proofs of Propositions \ref{propmotzkin} and
\ref{L1missing}.}
 For Motzkin
paths with arbitrary $t$,
the generating function $F(x)$ satisfies
$$F(x)=\frac{1}{1-tx-x^2F(x)}.$$
Applying Proposition \ref{l-quadraticTr} so $F_1=\mT(F)$ gives
$$\det( H_{n-1}(F_1))=\det (H_n(F)) \text{ where }
F_1(x)=\frac{1}{1-tx-x^2F_1(x)}.$$ Hence, $F(x)=F_1(x)$, and
consequently $\det (H_n(F(x)))=1$ for all $n$.

Whereas the Gessel-Viennot-Lindstr\"{o}m method leads  to a proof in the
shifted case for arbitrary $t$, as in
Proposition \ref{L1missing},  we have been able to use
the continued fractions technique only for $t=1$ and $t=2$.

For $t=1$  we will show that $(\det(H_n^1(F)))_{n\ge1} = (1,0,-1,-1,0,1,1,
\ldots)$, continuing with period 6.  Let $G_1(x) = (F(x)-1)/x$, so
that $\det(H_n^1(F))= \det(H_n(G_1)).$  Let $G_2 = \mT(G_1)$ and
 $G_3 =\mT(G_2)$, both under Proposition  \ref{l-quadraticTr}(ii).
 Since
 $$G_1(x) = \frac{ 1+x}{1-x-2x^2-x^3G_1(x)}$$
 with $d=0$, $k=3$, $u(x)=u_L(x)=1-2x$, $u_H = 0$, and $v(x) =-(1+x)^{-1}$,
 we find that
 $$G_2(x) =\frac{x}{1-x-2x^2-x^2(1+x)G_2(x)}$$
 with  $d=1$, $k=2$, $u(x)=u_L(x)=1-x-2x^2$, $u_H = 0$, and $v(x) =-(1+x)$.
 Applying Proposition \ref{l-quadraticTr}(ii)  shows
 \begin{eqnarray*}
 G_3(x) &=& -x^{-1}(-(1+x)) G_2(x)\\
 &=&  -x^{-1}(-(1+x)) (-x(-(1+x)^{-1})) G_1(x) \\
 &=& G_1(x)
 \end{eqnarray*}
 and $\det(H_{n-3}(G_3)) = -\det(H_{n-1}(G_2))=-\det(H_{n}(G_1))$,
 which yields the periodicity of the sequence of determinants.

 For $t=2$ we will show that $\det(H_n^1(F))  = n+1$ for $n \ge 1$.
 Define, $G_1$ to satisfy,
$$G_1(x) =\frac{2+x }{
       1- 2x - 2x^2
    - x^3 G_1(x) }.$$
One can easily see  that $G_1(x) =(F(x)-1)/x$  with $G_1(0) = u_1(0)^{-1} =
\det(H_1(G_1)) = \det(H_1^1(F)) = 2.$
For  $n \ge 2$,
 define, $G_n$ to satisfy,
$$  G_n(x)=\frac{(n-1)^2(n^2+n+x) }{
       (n^2-n)(n^2 - 2n^2x - 2x^2)
    - n^2(n^2-n+x)x^2 G_n(x) }.$$
By induction one can show that $G_{n} = \mT \circ \mT(G_{n-1})$ (under
Prop.~\ref{l-constantTr} then  under  Prop.~\ref{l-quadraticTr}),
and that $G_n(0) = u_n(0)^{-1} = (n-1)(n+1)/n^2.$
Also by induction and Proposition \ref{l-constantTr},
for $n\ge   2,$
\begin{eqnarray*}
\det(H_n(G_1))&=& \left[ 2^n\prod_{i=2}^{n-1}
\left(\frac{(i-1)(i+1)}{i^2}\right)^{n+1-i}\right] \det(H_1(G_n))\\
&=&  2^n\prod_{i=2}^{n}
\left(\frac{(i-1)(i+1)}{i^2}\right)^{n+1-i}
\end{eqnarray*}
which simplifies to $\det(H_n(G_1)) = n+1.$

\

{\bf Example 2:} {\sc Another proof of Proposition \ref{propschr}.}
 For large Schr\"oder paths  arbitrary $t$, we have
$$F(x)=\frac{1}{1-tx^2-x^2F(x)}.$$
Applying $\mT$ gives
$$\det (H_{n-1}(F_1))=\det (H_n(F)), \text{ where }
F_1(x)=\frac{1+t}{1+tx^2-x^2F_1(x)}.$$ Applying $\mT$ again, we
obtain
$$(1+t)^n\det (H_{n-1}(F_2))=\det (H_n(F_1)), \text{ where }
F_2(x)=\frac{1}{1-tx^2-x^2F_2(x)}.$$ This implies  $F_2=F$,
and hence the recurrence $\det (H_n (F)) = (1+t)^{n-1} \det
(H_{n-2}(F)),$ with initial condition $\det (H_1(F))=1$, and $\det
(H_2(F))=1+t$.

\

{\bf Example 3:} {\sc Another proof of Proposition \ref{proplargesch}}.
Consider the continued fraction
$$F(x) = \frac{1}{1-tx-xF(x)},$$
where $F(x)$ is the generating function for the Catalan numbers
for $t=0$ and the large Schr\"oder numbers for  $t=1$.

Under Proposition \ref{l-quadraticTr}(i) we have a unique $G_1$
such that
    $G_1(x) =F(x)$   and   $\det(H_{n-1}^1(G_1)) = \det(H_n(F)).$
Taking $G_2=(G_1(x)-1)/x=(F(x)-1)/x$, we have
$$ G_2(x) = \frac{(1+t)}{  1-(2+t)x -x^2G_2(x)} $$
where
$  \det(H_{n-1}(G_2)) = \det(H_{n-1}^1(F))$   and
     $u(x) =( 1-(2+t)x)/(1+t).$

Under Proposition \ref{l-constantTr}
 we have a unique $G_3$
$$ G_3(x) = \frac{1}{  1-(2+t)x -(1+t)x^2G_3(x)},$$
with $G_3(x) = G_2/(1+t) $ and
$\det(H_{n-1}(G_3)) = (1+t)^{-(n-1)}\det(H_{n-1}(G_2)).$

Under Proposition \ref{l-quadraticTr}(ii)  we have a unique
$G_4$ such that
$G_4(x) = (1+t)G_3(x)$ and
$\det(H_{n-2}(G_4)) = \det(H_{n-1}(G_3)).$

We see that $G_4(x) = G_2(x)$;
thus $\det(H_{n-1}(G_2))= (1+t)^{n-1}\det(H_{n-2}(G_2))$ with
$\det(H_{1}(G_2))= 1+   t$.  Hence $\det(H_{n}(F)) =
\det(H_{n-1}(G_2)) = (1+t)^{n(n-1)/2}.$

 \

{\bf Example 4:} {\sc Another proof of Proposition \ref{L2missing}.}
To compute  $\det(H_n^1(F))$, first we consider
$$H_n^1(F)=H_n(F_1), \text{ where } F_1=\frac{(t+1)x}{1-(2+t)x^2-x^3F_1}.$$
Applying $\mT$ shows that $\det(H_n(F_1))=-(1+t)^n \det( H_{n-2}(F_1)).$

\

{\bf Example 5:} For $\ell=3$, recall the functional equation
$$F_0(x)=\frac{1}{1-tx^3-x^2F_0(x)}.$$
For  arbitrary  $t$,
our transformation  gives more
   and more complicated expressions. This is not  surprising since the
   Hankel determinants do not factor nicely.
However, for $t=1$ and for $k=1,2,3$, the transformation
gives nice results similar to that of Proposition \ref{propperiod14}:
indeed,   sequences of $\det (H_n^k(F_0))$  also
 have period 14.  For $k=4$ there is an  interesting result.

\

{\bf Subexample 5i:}
The sequence for  $\det (H_n^1(F_0))$  starts with $ 0, -1, 0, 1, 1, 0, -1 ,
 0, 1, 0, -1, -1, 0, 1 $. If we define $F_1$ so that  $F_0(x)=1+xF_1(x)$,
then
$$
\det (H_n(F_1))  =\det
( H_n^1(F_0)), \text{ with } F_1=
\frac{x(x+1)}{1-2x^2-x^3-x^3F_1} \text{ and }d=1.
 $$
Then applying $\mT$
repeatedly so $\mT(F_i)=F_{i+1}$, we obtain

\begin{eqnarray*}
 \det(H_{n-2}(F_2)) &=& -\det (H_n(F_1)), \text{ where }
F_2=\frac{x}{(x+1)(1-x-x^2-x^3F_2)}\text{ and }d=1; \\
 \det(H_{n-2}(F_3))&=&-\det (H_n(F_2)), \text{ where }
F_3=\frac{1+x-x^2}{1-2x^2+x^3-x^3F_3}\text{ and }d=0; \\
 \det (H_{n-1}(F_4))&=&\det(H_{n}(F_3)), \text{ where }
F_3=\frac{x}{(1+x-x^2)(1-x-x^2F_3)}\text{ and }d=1 ;\\
  \det (H_{n-2}(F_5))&=&-\det(H_{n}(F_4)), \text{ where }
F_5= \frac{x(x+1)}{1-2x^2-x^3-x^3F_5}.
\end{eqnarray*}
The periodicity is established by  noticing that $F_5=F_1$ and $\det (H_{n-7}(F_5))=-\det (H_{n}(F_1)).$

\

{\bf Subexample 5ii:}
The sequence for  $\det (H_n^2(F_0))$  starts with $ 1,1,1,1,0,0,-1, -1,-1,-1,-1,0,0,1, $. If we define $G_0$ so that  $F_0(x)=1+x^2G_0(x)$,
then
$ \det (H_n(G_0))=\det(H_{n}^2(F_0))$,
$$G_0=\frac{1+x}{1-2x^2-x^3 - x^4G_0}.$$
One can establish the periodicity using Proposition \ref{l-quadraticTr}.
However,  this generating function has appeared in Transformation 3 of section \ref{proofperiod14}, where one can see that
\begin{equation}\label{beseen}
 \det (H_n(G_0)) =
-\det (H_{n+5}(F_0)).
\end{equation}

\

{\bf Subexample 5iii:}
The sequence for $\det (H_n^3(F_0))$  starts
with $1,-1,-1,0,0,0,-1,    -1,1,1,0,0,0,1$  and continues with period 14.
 The verification for this case uses Proposition \ref{l-quadraticTr}(ii) occasionally
 interspersed with
 Proposition \ref{l-constantTr}.
Here we will only sketch the verification. 
 By defining $F_1$ so that  $F_0(x)=1+x^2 +x^3F_1(x)$,
one finds that
$$F_1 =\frac{1+2x+x^2+x^3}{1-2x^2- x^3 -2x^4 -x^5 F_1}.$$
For the first transformation, with $F_2 = \mT F_1$, we find
 $$F_2 = \frac{1-2x+x^3}{ -1 + 4x^2 +x^3 + 2x^4 - x^2(1+2x+x^2+x^3)F_2},$$
in which $u(x)= (-1 + 4x^2 +x^3 + 2x^4)/(1-2x+x^3) $.  Now, since
 $ u(0) = -1$,  one needs to apply Proposition \ref{l-constantTr} for the
next transformation. One proceeds until a generating function equal to
  $F_1$ appears to establish the periodicity.  We remark that $d=0$ for each transformation until the
 final one which uses
Proposition \ref{l-quadraticTr}(ii) with $d=3$
(This corresponds to  a fourth order block).

\

{\bf Subexample 5iv:}
The sequence for $\det (H_n^4(F_0))$  begins with $$2, 3, 4, 0, 0, -4, -5, -6 ,-7,
 -8, 0, 0, 8, 9, 10, 11, 12, 0, 0, -12, -13, -14, -15, -16, 0,
     0,16,\dots.$$
For $n \ge 8$, an essence of periodicity can be gleaned from the recurrence
$$\det(H_{n}^4(F_0))=4\det ( H_{n-1}(F_0))- \det (H_{n-7}^4(F_0)),$$
for which we sketch a proof, often omitting the functional equations.

We will be applying the transformation $\mT$ eight times,
alternating its definition to be first under
Proposition \ref{l-constantTr} and then under Proposition \ref{l-quadraticTr}(ii).
Let $F_1$ satisfy
$F_0=1+x^2+x^3+x^4F_1$.  Hence,  $\det (H_n(F_1))=\det ( H_n^4(F_0))$, and
$$F_1 = \frac{2+3x+2x^2+2x^3+x^4}{1-2x^2-x^3-2x^4-2x^5-x^6F_1}.$$

 Here $u(0)=\frac{1}{2}$, where $u(x)$ is for $F_1$.  Thus, with $F_2=\mT F_1$,
 $\det( H_n(F_2))= (\frac{1}{2})^n\det (H_n(F_1))$.
Now $d=0$, where $d$ is for $F_2$.  With $F_3=\mT F_2$,
 $\det (H_{n-1}(F_3))= \det (H_n(F_2))$.

Here $u(0)=\frac{4}{3}$, where $u(x)$ is for $F_3$.  Thus, with $F_4=\mT F_3$,
 $\det ( H_{n-1}(F_4))= (\frac{4}{3})^{n-1}\det  (H_{n-1}(F_3))$.
Now $d=0$, where $d$ is for $F_4$.  With $F_5=\mT F_4$,
 $\det  (H_{n-2}(F_5))= \det  (H_{n-1}(F_4))$.

Here $u(0)=\frac{9}{8}$, where $u(x)$ is for $F_5$.  Thus, with $F_6=\mT F_5$,
 $\det  (H_{n-2}(F_6) )= (\frac{9}{8})^{n-2}\det  ( H_{n-2}(F_5))$.
Now $d=0$, where $d$ is for $F_6$.  With $F_7=\mT F_6$,
 $\det (H_{n-3}(F_7))= \det (H_{n-2}(F_6))$.

Here $u(0)=\frac{4}{3}$, where $u(x)$ is for $F_7$.    Thus, with $F_8=\mT F_7$,
 $\det (H_{n-3}(F_8))= (\frac{4}{3})^{n-3}\det (H_{n-3}(F_7))$.

Now $d=2$, where $d$ is for $F_8$.  With $F_9=\mT F_8$,
 $\det( H_{n-6}(F_9))= - \det( (H_{n-3}(F_8))  =  \\
 (\frac{5}{4},\frac{6}{4},\frac{7}{4},\frac{8}{4},0,0,-\frac{8}{4},-\frac{9}{4},-\frac{10}{4},\ldots).$

Thus, (surprisingly)
\begin{eqnarray}
\det (H_{n-6}(F_9)) & =& -(\frac{1}{2})^{n}(\frac{4}{3})^{n-1}
(\frac{9}{8})^{n-2}
(\frac{4}{3})^{n-3}
 \det (H_{n}^4(F_0)) \notag
\\
 &=& -\frac{1}{4} \det ( H_{n}^4(F_0)) \label{ident1}
\end{eqnarray}
Moreover,$$ F_9=\frac{20+16x-8x^2-4x^3+x^4}{8(2-4x^2-2x^3+x^4)-16x^4F_9}=
\frac{5}{4}+x+2x^2+3x^3+6x^4+10x^5+\cdots$$
It is easily verified that $F_9(x)$ and  $\frac{1}{4}+G_0(x)$, where $G_0$ appears in Subexample 5ii,   satisfy the same functional equation, and hence are equal.
Therefore,
\begin{eqnarray*}
\det (H_{n-6}(F_9)) &=& \left[\frac{xF_9(x)-yF_9(y)}{x-y} \right]_{n-6} \\
&=& \left[\frac{1}{4}+ \frac{xG_0(x)-yG_0(y)}{x-y} \right]_{n-6}    \\
&=& \frac{1}{4}\det (H_{n-7}^4(F_0)) +\det ( H_{n-6}(G_0)) \\
\end{eqnarray*}
where $\frac{1}{4}\det( H_{n-7}^4(F_0)) $ is  $\frac{1}{4}$ times the determinant of  the 1,1-minor of $H_{n-6}(G_0)$, equivalently of $H_{n-6}^2(F_0)$.
 Combining this with identity (\ref{ident1})
 and noting $\det(H_n(G_0)) =
-\det (H_{n+5}(F_0))$ from (\ref{beseen})
 proves the initial recurrence of this subexample.

\

{\bf Acknowledgments.} We, particularly Xin,  are most grateful  to
Ira Gessel.
We also appreciate the suggestions made
by Christian Krattenthaler, Lou Shapiro, and referees.

\bibliographystyle{amsalpha}

\begin{thebibliography}{A}

\bibitem{Aigner}
M.~Aigner, Motzkin Numbers, {\it Europ. J. Comb.} 19
(1998), 663-675.

\bibitem{BanderierSchwer}  C.~Banderier and  S.~Schwer,
 Why Delannoy numbers?, {\it J. Stat. Plan. and Infer.},  135 (2005), 40-54.


\bibitem{Bressoud} D.~Bressoud, {\it Proofs and Confirmations: The Story of the Alternating Sign Matrix Conjecture},
Cambridge University Press, 1999.


\bibitem{BrualdiKirkland}
R.~A.~Brualdi and S.~Kirkland,  Aztec diamonds and digraphs, and Hankel
determinants of Schr\"oder numbers, {\it J. Comb. Th., ser. B}, 94 (2005), 334-351.

\bibitem{GV}
I.~M.~Gessel and G.~X.~Viennot, Binomial determinants, paths, and hook
length formulae. {\it Adv. in Math.}  58 (1985), no. 3, 300-321.

\bibitem{xingessel2}
I.~M.~Gessel and G.~Xin, \emph{The generating function of ternary
trees and
  continued fractions}, Electron. J. Combin., to appear,\\
ArXiv:math.CO/0505217.

\bibitem{cfraction}
W.~B.~Jones and W.~J.~Thron, \emph{Continued Fractions: Analytic
Theory and Applications},
  Encyclopedia of Mathematics and its Applications, vol.~11, Addison-Wesley
  Publishing Co., Reading, Mass., 1980.

\bibitem{kratt}
 C.~Krattenthaler,
 Advanced determinant calculus: a complement, {\it Linear Algebra
     Appl.} 411 (2005), 68-166

The Andrews Festschrift (Maratea, 1998). {\it Sem. Lothar.
Combin.}  42 (1999), Art. B42q, 67 pp. (electronic).

\bibitem{Pergolaetal} E.~Pergola, R.~Pinzani,   S.~Rinaldi, and R.~A.~Sulanke,
Lattice Path Moments by Cut and Paste. {\it
Adv. in Appl. Math.} 30 (2003), no. 1-2, 208-218.

\bibitem{Shapiro}  Shapiro, L. W. A Catalan triangle. {\it Discrete Math.} 14 no. 1  (1976),  83--90.

\bibitem{sulMomentsgenMot}  R.~A.~Sulanke,
Moments of generalized Motzkin paths,
   {\it J. Int. Seq.},  3 (2000), Article 00.1.1

\bibitem{SulDel}
R.~A.~Sulanke,
Objects counted by the central Delannoy numbers. {\it J. Int.
Seq.} 6 (2003), no. 1, Art. 03.1.5, 19 pp.

\bibitem{Viennot} G.~Viennot,
  Une th\'eorie combinatoire des polyn\^omes orthogonaux g\'en\'eraux,
  Lecture notes, Univ. Quebec, Montreal, Que., 1983.

\bibitem{wall}
H.~S.~Wall, \emph{Analytic Theory of Continued Fractions}, Van
Nostrand, New York, 1948.


\end{thebibliography}

\end{document}